\documentclass[12pt,reqno]{amsart}
\usepackage{amsmath, amsthm, amssymb}
\usepackage{color}

\numberwithin{equation}{section}

\newcommand{\mmod}[1]{\,\,(\text{\rm{mod}}\,\,#1)}

\def\bfa{{\boldsymbol a}}
\def\bfb{{\boldsymbol b}}

\def\phi{{\varphi}}

\def\le{\leqslant} \def\ge{\geqslant}

\newtheorem{thm}{Theorem}[section]

\newtheorem{conj}[thm]{Conjecture}

\newtheorem{lem}[thm]{Lemma}
\newtheorem{coro}[thm]{Corollary}

\def\pf{\noindent{\it Proof.} }
\def\qed{\nopagebreak\hfill{\rule{4pt}{7pt}}
\medbreak}

\numberwithin{equation}{section}

\def\qed{\nopagebreak\hfill{\rule{4pt}{7pt}}
\medbreak}

\newlength{\boxedparwidth}
  {\begin{center} \begin{tabular}{|@{\hspace{.315in}}c@{\hspace{.15in}}|}
                  \hline \\ \begin{minipage}[t]{\boxedparwidth}
                  \setlength{\parindent}{.25in}}%
  {\end{minipage} \\ \\ \hline \end{tabular} \end{center}}
\allowdisplaybreaks
\begin{document}

\title[Generalized Lambert Series Identities]
{Generalized Lambert Series Identities and Applications in Rank Differences}
\author[Bin Wei]{Bin Wei}
\address{Center for Applied Mathematics, Tianjin University, Tianjin 300072, P.R. China}
\email{bwei@tju.edu.cn}
\author[Helen W.J. Zhang]{Helen W.J. Zhang}
\address{Center for Applied Mathematics, Tianjin University, Tianjin 300072, P.R. China}
\email{wenjingzhang@tju.edu.cn}
\thanks{The authors are supported by NSFC (Grant NO. 11701412).}

\maketitle

\noindent {\bf Abstract.}
In this article, we prove two identities of generalized Lambert series.
By introducing what we call $\mathcal{S}$-series, we establish relationships between multiple generalized Lambert series and multiple infinite products.
Compared with Chan's work, these new identities are useful in generating various formulas for generalized Lambert series with the same poles.
Using these formulas, we study the 3-dissection properties of ranks for overpartitions modulo 6.
In this case, $-1$ appears as a unit root, so that double poles occur.
We also relate these ranks to the third order mock theta functions $\omega(q)$ and $\rho(q)$.

\setcounter{section}{-1}
\section{Notations}

Throughout this article, we use the common $q$-series notations associated with infinite products:
\begin{align*}
&(a)_\infty:=(a;q)_\infty:=\prod_{n=0}^\infty(1-aq^n),
&&(a_1,a_2,\ldots,a_k)_\infty:=(a_1)_\infty\cdots(a_k)_\infty,
\\[5pt]
&[a]_\infty:=(a,q/a)_\infty,
&&[a_1,a_2,\ldots,a_k]_\infty:=[a_1]_\infty\cdots[a_k]_\infty,
\\[5pt]
&j(z;q):=(z;q)_\infty (q/z;q)_\infty (q;q)_\infty,
&&J_{a,m}:=j(q^a;q^m),\ \  J_m:=(q^m;q^m)_\infty.
\end{align*}
For the sake of convergence, we always assume that $|q|<1$. Also, we adopt a notation due to D. B. Sears \cite{Sears-1951}:
\begin{align*}
&F(b_1,b_2,\ldots,b_m)+\mathrm{idem}(b_1;b_2,\ldots,b_m)
\\
&:=F(b_1,b_2,\ldots,b_m)+F(b_2,b_1,b_3,\ldots,b_m)
+\cdots+F(b_m,b_2,\ldots,b_{m-1},b_1).
\end{align*}

\section{Introduction}

A Lambert series, named for Johann Heinrich Lambert, takes the form
$$
\sum_{n=1}^\infty a_n\frac{q^n}{1-q^n},
$$
where $\{a_n\}$ is any set of real or complex numbers. A generalized Lambert series allows more general exponents in both numerators and denominators.
Such series are often useful in obtaining formulas for various generating functions, since the denominators can be expanded as a geometric series.
Expanded generalized Lambert series are naturally linked with infinite products. For example, Chan \cite{Chan-2005} proved three generalized Lambert series expansions for infinite products.
One of the theorems concerning $r+1$ poles in generalized Lambert series is stated as following.
\begin{lem}\label{Chan-Thm-2.2}
For non-negative integers $r<s$, we have
\begin{align*}
&\frac{(a_1q,q/a_1,\ldots,a_rq,q/a_r,q,q)_\infty}
{[b_1,b_2,\cdots,b_s]_\infty}
\\
&\indent=
\frac{[a_1/b_1,\cdots,a_r/b_1]_\infty}
{[b_2/b_1,\cdots,b_s/b_1]_\infty}
\sum_{n=-\infty}^\infty\frac{(-1)^{(s-r)n+r}q^{(s-r)n(n+1)/2}b_1^ra_1^{-1}\cdots a_r^{-1}}
{(1-b_1q^n)(1-b_1q^n/a_1)\cdots(1-b_1q^n/a_r)}
\\
&\indent\quad \times\left(\frac{a_1\cdots a_rb_1^{s-r-1}q^r}
{b_2\cdots b_s}\right)^n
+\mathrm{idem}(b_1;b_2,\ldots,b_s).
\end{align*}
For $r=s$, this is true provided that $|q|<|\frac{a_1\cdots a_r}{b_1\cdots b_s}|<|q^{-r}|$.
\end{lem}

Using these theorems, Chan provided brief proofs for amounts of beautiful and useful identities.
Particularly, when taking $s=3$ and $r=0$, Lemma \ref{Chan-Thm-2.2} delivers the key identity used by Atkin and Swinnerton-Dyer \cite{Atkin-Swinnerton-Dyer-1954} in proving Ramanujan's famous partition congruences.

One limitation of applying Chan's theorems is that, the exponents in numerators are partially or totally determined by the poles.
In this article, we prove two other generalized Lambert series identities.
First, for a sequence $\bfa=(a_1,\ldots,a_r)$, we define series $\mathcal{S}(a_1,\ldots,a_r)$ as
\begin{equation}\label{F(a,q)}
\mathcal{S}(a_1,\ldots,a_r):=\mathcal{S}(a_1,\ldots,a_r;q)=
\sum_{u=1}^r\sum_{n=0}^\infty\left(\frac{1}{1-a_uq^n}
-\frac{1}{1-a_u^{-1}q^{n+1}}\right).
\end{equation}
We also write $\mathcal{S}(\bfa)=\mathcal{S}(a_1,\ldots,a_r)$ for brevity.
The following theorem concern identity of generalized Lambert series with double poles.

\begin{thm}\label{main-result}
Let $\bfa=(a_1,\ldots,a_r)$ and $\bfb=(b_1,\ldots,b_s)$. Then for non-negative integers $r<s$, we have
\begin{align}\label{R}
&\frac{(q)_\infty^2[a_1,\cdots,a_r]_\infty}
{[b_1,\cdots,b_s]_\infty}\left(1-\mathcal{S}(\bfa)+\mathcal{S}(\bfb)\right)
\nonumber\\[5pt]
&\quad=\frac{[a_1/b_1,\cdots,a_r/b_1]_\infty}
{[b_2/b_1,\cdots,b_s/b_1]_\infty}\sum_{n=-\infty}^\infty
\frac{(-1)^{(s-r)n}q^{(s-r)n(n+1)/2}}{(1-b_1q^n)^2}
\left(\frac{a_1\cdots a_rb_1^{s-r-1}}{b_2\cdots b_s}\right)^n
\nonumber\\[5pt]
&\indent\indent\quad+\mathrm{idem}(b_1;b_2,\cdots,b_s).
\end{align}
For $r=s$, this is true provided that $|q|<|\frac{a_1\cdots a_r}{b_1\cdots b_s}|<1$.
\end{thm}

A similar identity concerning generalized Lambert series with single poles is also given in \S2.
Theorem \ref{main-result} is aimed at decoupling parameters in $\bfa$ from the denominators.
Therefore, it is helpful in generating various identities concerning generalized Lambert series with the same poles.

The generalized Lambert series $\mathcal{S}$ defined in (\ref{F(a,q)}) appears as an encumbrance in our identities for infinite products.
Though, we provide an algorithm to show that $\mathcal{S}(\pm q^m;q^n)$ with $m$, $n$ integers can be expanded as sums of infinite products.
Therefore, our main result Theorem \ref{main-result} establishes a bridge between multiple infinite products and multiple generalized Lambert series.
For example, we show that the following identity holds in \S4:
$$
\nonumber\sum_{n=-\infty}^\infty\frac{(-1)^n q^{3n^2+3n}}{(1+q^{3n+1})^2}
+\sum_{n=-\infty}^\infty\frac{(-1)^n q^{3n^2+3n}}{(1-q^{3n+1})^2}
=\frac{4J^3_{6}}{3J_{2}}+\frac{J_{3,6}^2J_{6}^6}{2J_{1,6}^2J_2^2}+\frac{J_{1,6}^6J_{2}^2J_{3,6}^2}{6J_{6}^6}.
$$

The motivation of establishing new identities for generalized Lambert series arises in the study of the series
$$
\overline{R}(-1;q)=\frac{4(-q)_\infty}{(q)_\infty}\sum_{n=-\infty}^\infty
\frac{(-1)^nq^{n^2+n}}{(1+q^{n})^2}.
$$
Bringmann and Lovejoy \cite{Bringmann-Lovejoy-2007} proved that $\overline{R}(-1;q)$ is the holomorphic part of a harmonic weak Maass form of weight $3/2$.
They also pointed out that this is the most complicated case among $\overline{R}(z;q)$ since double poles occur.
In this article, we use Theorem \ref{main-result} to give the 3-dissection properties of $\overline{R}(-1;q)$.

Recall that an overpartition of positive integer $n$, denoted by $\overline{p}(n)$, is a partition of $n$ where the first occurrence of each distinct part may be overlined.
Particularly, we set $\overline{p}(0)= 1$.
The rank of an overpartition was introduced by Lovejoy \cite{Lovejoy-2005} as the largest part minus the number of parts.
Let $\overline{N}(m,n)$ denote the number of overpartitions of $n$ with the rank $m$, and let $\overline{N}(s,\ell,n)$ denote the number of overpartitions of $n$ of rank congruent to $s$ modulo $\ell$.
Lovejoy gave a generating function of $\overline{N}(m,n)$
\begin{align}\label{GF}
\overline{R}(z;q)
&:=\sum_{n=0}^\infty\sum_{m=-\infty}^{\infty}\overline{N}(m,n)z^mq^n
\nonumber\\[5pt]
&=\frac{(-q)_\infty}{(q)_\infty}\left\{1+2\sum_{n=1}^\infty
\frac{(1-z)(1-z^{-1})(-1)^nq^{n^2+n}}{(1-zq^n)(1-z^{-1}q^n)}\right\}.
\end{align}

Rank differences between different residues are widely studied, where identities of generalized Lambert series usually play key roles.
In \cite{Lovejoy-Osburn-2008}, Lovejoy and Osburn gave formulas for the full rank differences $\overline{N}(s, \ell, \ell n + d)-\overline{N}(t, \ell, \ell n + d)$ for $\ell=3,5$, in terms of infinite products and generalized Lambert series.
The modulus $7$ have been determined by Jennings-Shaffer \cite{Jennings-Shaffer-2016}.
Besides, when considering even moduli, only special linear combinations of rank differences can be obtained previously.
In \cite{Ji-Zhang-Zhao-2017}, Ji, Zhang and Zhao studied 3-dissection properties of the form
\begin{equation}\label{Ji-Zhang-Zhao}
\sum_{n=0}^{\infty}(\overline{N}(0,6,n)+\overline{N}(1,6,n)-\overline{N}(2,6,n)-\overline{N}(3,6,n))q^n.
\end{equation}
The difficulty of providing full rank differences lies in the truth that $-1$ is a unit root of even moduli, so that $\overline{R}(-1;q)$ arises naturally.
Similar situation happens in related problems, which are associated with various types of ranks (such as crank, $M_2$-rank, etc.) for different types of partitions (see \cite{Andrews-Lewis-2000,Garvan-1988,Garvan-1990,Mao-2013,Mao-2015} for example).

As a consequence of successfully handling double poles, we are now able to handle the full rank differences associated with even moduli.
In fact, we gave the formulas for each residue, i.e., each term in \eqref{Ji-Zhang-Zhao} instead.
Here we take the 3-dissection properties of ranks of overpartitions modulo $6$ as an example.
Let
\begin{equation}\label{r_s(d)}
\overline{r}_s(d)=\sum_{n=0}^\infty\overline{N}(s,6,3n+d)q^n.
\end{equation}
When $d=2$, we have the following theorem.
\begin{thm}\label{rank-diff-2}
We have
\begin{align*}
\overline{r}_0(2)
&=\frac{2J_{6}^{12}}{3J_{1,6}^6J_{2}^4J_{3,6}^3}
-\frac{4J^3_{6}}{3J_{2}J_{3,6}}
+\frac{2J_{1,6}^6J_2^2J_{3,6}}{3J_{6}^6}
\nonumber\\[5pt]
&\quad\quad\quad
-\frac{4}{J_{3,6}}
\sum_{n=-\infty}^\infty
\frac{(-1)^nq^{3n^2+3n}}{(1+q^{3n+1})^2}
+\frac{4}{J_{3,6}}\sum_{n=-\infty}^\infty\frac{(-1)^n q^{3n^2+3n}}{1+q^{3n+1}},\\
\overline{r}_1(2)
&=\frac{2J_{6}^{12}}{3J_{1,6}^6J_{2}^4J_{3,6}^3}
+\frac{2J_{6}^3}{J_{2}J_{3,6}}
-\frac{2J_{1,6}^6J_2^2J_{3,6}}{3J_{6}^6}
\nonumber\\[5pt]
&\quad\quad\quad+\frac{4}{J_{3,6}}
\sum_{n=-\infty}^\infty
\frac{(-1)^nq^{3n^2+3n}}{(1+q^{3n+1})^2}
-\frac{4}{J_{3,6}}\sum_{n=-\infty}^\infty\frac{(-1)^n q^{3n^2+3n}}{1+q^{3n+1}},\\
\overline{r}_2(2)
&=\frac{2J_{6}^{12}}{3J_{1,6}^6J_{2}^4J_{3,6}^3}
+\frac{2J^3_{6}}{3J_{2}J_{3,6}}
+\frac{2J_{1,6}^6J_2^2J_{3,6}}{3J_{6}^6}
\nonumber\\[5pt]
&\quad\quad\quad
-\frac{4}{J_{3,6}}
\sum_{n=-\infty}^\infty
\frac{(-1)^nq^{3n^2+3n}}{(1+q^{3n+1})^2}
+\frac{2}{J_{3,6}}\sum_{n=-\infty}^\infty\frac{(-1)^n q^{3n^2+3n}}{1+q^{3n+1}},\\
\overline{r}_3(2)
&=\frac{2J_{6}^{12}}{3J_{1,6}^6J_{2}^4J_{3,6}^3}
-\frac{4J_{6}^3}{J_2J_{3,6}}
-\frac{2J_{1,6}^6J_2^2J_{3,6}}{3J_{6}^6}
\nonumber\\[5pt]
&\quad\quad\quad
+\frac{4}{J_{3,6}}\sum_{n=-\infty}^\infty
\frac{(-1)^nq^{3n^2+3n}}{(1+q^{3n+1})^2}.
\end{align*}
\end{thm}

The formulas for residues $d=0$ and $1$ are listed in \S5.
One should not be surprised of simultaneous occurrences of terms containing denominators $(1+q^{3n+1})$ and $(1+q^{3n+1})^2$, since double poles exist.
These explicit formulas suggest inequalities of ranks between different residues, such as $\overline{N}(1,6,3n+2)\ge\overline{N}(3,6,3n+2)$.
A conjecture on total ordering will also be discussed in \S5.

In \cite{Hickerson-Mortenson-2014}, Hickerson and Mortenson showed that a mock theta function can be expressed in terms of Appel-Lerch sums.
Inspired by their work, we establish a relation between the third order mock theta functions $\omega(q)$ and $\rho(q)$ and  the ranks of overpartitions modulo 6, where $\omega(q)$ and $\rho(q)$ are defined by
{\rm \cite{Watson-1936}:}
$$\omega(q)=\sum_{n=0}^\infty\frac{q^{2n(n+1)}}{(q;q^2)_{n+1}^2} \quad \text{and} \quad
\rho(q)=\sum_{n=0}^\infty\frac{q^{2n(n+1)}(q;q^2)_{n+1}}{(q^3;q^6)_{n+1}}.$$

\begin{thm}\label{mock}
We have
\begin{align}
\overline{r}_0(2)+\overline{r}_3(2)
&=\frac{4}{9}\rho(q)-\frac{16}{9}\omega(q)+M(q),
\\[3pt]
\overline{r}_1(2)-\overline{r}_3(2)
&=2\omega(q),
\\[3pt]
\overline{r}_2(2)+\overline{r}_3(2)
&=-\frac{2}{9}\rho(q)-\frac{10}{9}\omega(q)+M(q),
\end{align}
where $M(q)$ is (explicit) weakly holomorphic modular form given by:
\[M(q)=\frac{4J_{6}^{12}}{3J_{1,6}^6J_{2}^4J_{3,6}^3}.\]
\end{thm}

This paper is organized as follows.
In \S 2, we derive the main theorems by discussing poles in Chan's identities.
In \S 3, we introduce an algorithm for $\mathcal{S}$-series, which helps transform $\mathcal{S}$-series into sums of infinite products.
In \S 4, we use our new identities to generate some formulas concerning the 3-dissections of generalized Lambert series.
These formulas help establish 3-dissection properties of ranks for overpartitions modulo 6, in \S 5.
Finally, we prove the relations between the ranks of overpartitions and mock theta functions in \S 6.

\section{Proofs of Main Theorems}

We start with the following lemma, where we have made slight variants in the subscripts of parameters.
\begin{lem}[Chan \cite{Chan-2005}]\label{Chan-Thm-3.2}
For non-negative integers $r<s$, we have
\begin{align*}
&\frac{(a_0q,q/a_0,q,q)_\infty[a_1,a_2,\cdots,a_r]_\infty}
{[b_0,b_1,\cdots,b_s]_\infty}
\\
&\indent\indent=
\frac{[a_0/b_0,a_1/b_0,\cdots,a_r/b_0]_\infty}
{[b_1/b_0,b_2/b_0,\cdots,b_s/b_0]_\infty}
\sum_{n=-\infty}^\infty\frac{(-1)^{(s-r)n+1}q^{(s-r)n(n+1)/2}b_0a_0^{-1}}
{(1-b_0q^n)(1-b_0q^n/a_0)}
\\
&\indent\indent\quad \times\left(\frac{a_0a_1\cdots a_rb_0^{s-r-1}q}
{b_1\cdots b_s}\right)^n
+\mathrm{idem}(b_0;b_1,b_2,\ldots,b_s).
\end{align*}
For $r=s$, this is true provided that $|q|<|\frac{a_1\cdots a_r}{b_1\cdots b_s}|<1$.
\end{lem}

Compared with Lemma \ref{Chan-Thm-3.2}, Theorem \ref{main-result} is aimed at decoupling parameters in $\bfa$ from poles in generalized Lambert series, making it convenient to control orders of $q$ in numerators.
When generating identities in special forms, this also permits us to save a parameter in $\bfb$, and so a term of generalized Lambert series.\\

\noindent\emph{Proof of Theorem 1.2.} Briefly speaking, Theorem \ref{main-result} follows from setting $a_0=1$ and $b_0=q$ in Lemma \ref{Chan-Thm-3.2}.
Obviously this would results in double poles in both sides, so we need to compute the limits at $b_0=q$.
First replacing $a_0$ by $1$ in Lemma \ref{Chan-Thm-3.2}, we have
\begin{align}\label{Chan-gene}
&\frac{(q)_\infty^4[a_1,\cdots,a_r]_\infty}{[b_0,b_1,\cdots,b_s]_\infty}
=\frac{[b_0^{-1},a_1/b_0,\cdots,a_r/b_0]_\infty}
{[b_1/b_0,\cdots,b_s/b_0]_\infty}
\nonumber\\
&\quad\quad\quad\times\sum_{n=-\infty}^\infty\frac{(-1)^{(s-r)n+1}
q^{(s-r)n(n+1)/2}b_0}{(1-b_0q^n)^2}
\left(\frac{a_1\cdots a_rb_0^{s-r-1}q}{b_1\cdots b_s}\right)^n
\nonumber\\
&\quad\quad\quad+\mathrm{idem}(b_0;b_1,\cdots,b_s).
\end{align}

Denote the term on the left-hand side of (\ref{Chan-gene}) by $L$ and those on the right-hand side by $R_0,\ldots,R_s$ respectively, which is
\[L=R_0+R_1+\cdots+R_s.\]
For the right-hand side, the pole $b_0=q$ occurs only in the term $R_0$. So we may set $b_0\rightarrow q$ directly in other terms. As for $R_1$, we have
\begin{align}\label{R_1}
\nonumber\lim_{b_0\rightarrow q}R_1
&=\frac{[a_1/b_1,\cdots,a_r/b_1]_\infty}
{[b_2/b_1,\cdots,b_s/b_1]_\infty}\\
&\quad\times\sum_{n=-\infty}^\infty\frac{(-1)^{(s-r)n}q^{(s-r)n(n+1)/2}}{(1-b_1q^n)^2}
\left(\frac{a_1\cdots a_rb_1^{s-r-1}}{b_2\cdots b_s}\right)^n,
\end{align}
which is the first term of the right-hand side in (\ref{R}).
Thus it remains to show that
\begin{align}\label{L-R_0}
\lim_{b_0\rightarrow q}(L-R_0)=\frac{(q)_\infty^2[a_1,\cdots,a_r]_\infty}
{[b_1,\cdots,b_s]_\infty}\left(1-\mathcal{S}(\bfa)+\mathcal{S}(\bfb)\right).
\end{align}

We separate terms containing poles from $L$ and $R_0$ successively. We begin with rewriting $L$ and $R$ as
\begin{align*}
L=\frac{(q)_\infty^4[a_1,\cdots,a_r]_\infty}
{(b_0,b_0^{-1}q^2)_\infty[b_1,\cdots,b_s]_\infty}
\cdot\frac{b_0}{b_0-q},
\end{align*}
and
\begin{align}
R_0\label{R_0}
\nonumber&=\frac{(1-b_0^{-1}q)(b_0,b_0^{-1}q^2)_\infty
[a_1q/b_0,\cdots,a_rq/b_0]_\infty}
{[b_1q/b_0,\cdots,b_sq/b_0]_\infty}
\\[6pt]
&\quad\times\sum_{n=-\infty}^\infty
\frac{(-1)^{(s-r)(n+1)}q^{(s-r)n(n+1)/2}b_0}{(1-b_0q^n)^2}
\left(\frac{a_1\cdots a_rb_0^{s-r}}{b_1\cdots b_s}\right)^{n+1}\left(\frac{q}{b_0}\right)^n.
\end{align}
It is easy to see that, in the generalized Lambert series in $R_0$, poles occurs only when $n=-1$.
Considering the factor $(1-b_0^{-1}q)$, other terms vanish when setting $b_0\rightarrow q$. Thus, in (\ref{L-R_0}), we have
\begin{align}\label{L^*-R_0^*}
\nonumber\lim_{b_0\rightarrow q}(L-R_0)&=\lim_{b_0\rightarrow q}\frac{1}{b_0-q}
\left(\frac{(q)_\infty^4[a_1,\cdots,a_r]_\infty}
{(b_0,b_0^{-1}q^2)_\infty[b_1,\cdots,b_s]_\infty}
\cdot b_0\right.
\\[6pt]
&\nonumber\quad\quad\quad\left.-\frac{(b_0,b_0^{-1}q^2)_\infty
[a_1q/b_0,\cdots,a_rq/b_0]_\infty}
{[b_1q/b_0,\cdots,b_sq/b_0]_\infty}\cdot q\right)
\\[6pt]
&\nonumber=\lim_{b_0\rightarrow q}\frac{{\rm d}}{{{\rm d}} b_0}
\left(\frac{(q)_\infty^4[a_1,\cdots,a_r]_\infty}
{(b_0,b_0^{-1}q^2)_\infty[b_1,\cdots,b_s]_\infty}
b_0\right.
\\[6pt]
&\nonumber\quad\quad\quad\left.-\frac{(b_0,b_0^{-1}q^2)_\infty
[a_1q/b_0,\cdots,a_rq/b_0]_\infty}
{[b_1q/b_0,\cdots,b_sq/b_0]_\infty}q\right)
\\[6pt]
&:=\lim_{b_0\rightarrow q}\frac{{\rm d}}{{{\rm d}} b_0}(L^*-R_0^*),
\end{align}
where the penultimate equation follows by L'H\^{o}pital's rule.

For $L^{*}$, It is easy to obtain
\begin{align}\label{L^{*}}
\lim_{b_0\rightarrow q}\frac{{\rm d}L^{*}}{{\rm d}b_0}&=\frac{(q)_\infty^2[a_1,\cdots,a_r]_\infty}{[b_1,\cdots,b_s]_\infty}.
\end{align}

For $R_0^*$, we have
\begin{align*}
\lim_{b_0\rightarrow q}R_0^*=\frac{q(q)_\infty^2[a_1,\cdots,a_r]_\infty}
{[b_1,\cdots,b_s]_\infty}.
\end{align*}
It follows by taking the logarithmic derivative that
\begin{align*}
\lim_{b_0\rightarrow q}\frac{{\rm d}\log R_0^*}{{\rm d}b_0}=\frac{\mathcal{S}(\bfa)-\mathcal{S}(\bfb)}{q},
\end{align*}
where $\mathcal{S}$ is defined in \eqref{F(a,q)}.
Therefore,
\begin{align}\label{R_0^*}
\lim_{b_0\rightarrow q}\frac{{\rm d}R_0^*}{{\rm d}b_0}
&=\lim_{b_0\rightarrow q}\left(R_0^*\frac{{\rm d}\log R_0^*}{{\rm d}b_0}\right)
\nonumber\\[5pt]
&=\frac{(q)_\infty^2[a_1,\cdots,a_r]_\infty}
{[b_1,\cdots,b_s]_\infty}\left(\mathcal{S}(\bfa)-\mathcal{S}(\bfb)\right).
\end{align}
Thus we complete the proof by substituting (\ref{L^{*}}) and (\ref{R_0^*}) into (\ref{L^*-R_0^*}).
\qed

Chan \cite{Chan-2005} also proved the following identity concerning generalized Lambert series with single poles.

\begin{lem}\label{Chan-Thm-2.1}
For non-negative integers $r<s$, we have
\begin{align}\label{Chan-2.1}
\frac{[a_1,\cdots,a_r]_\infty(q)_\infty^2}
{[b_0,b_1,\cdots,b_s]_\infty}
&=\frac{[a_1/b_0,\cdots,a_r/b_0]_\infty}
{[b_1/b_0,\cdots,b_s/b_0]_\infty}\nonumber\\
&\times\sum_{n=-\infty}^\infty\frac{(-1)^{(s-r+1)n}q^{(s-r+1)n(n+1)/2}}{1-b_0q^n}
\left(\frac{a_1\cdots a_rb_0^{s-r}}
{b_1\cdots b_s}\right)^n
\nonumber\\
&\indent\indent\quad
+\mathrm{idem}(b_0;b_1,\ldots,b_s).
\end{align}
For $r=s$, this is true provided that $|q|<|\frac{a_1\cdots a_r}{b_1\cdots b_s}|<1$.
\end{lem}

Similarly by taking $b_0\rightarrow q$, we obtain the following theorem.

\begin{thm}\label{single-pole}
Let $\bfa=(a_1,\ldots,a_r)$ and $\bfb=(b_1,\ldots,b_s)$.  Then for non-negative integers $r\le s$, we have
\begin{align}\label{R'}
&\frac{[a_1,\cdots,a_r]_\infty}
{[b_1,\cdots,b_s]_\infty}\left(1-\mathcal{S}(\bfa)+\mathcal{S}(\bfb)\right)
\nonumber\\[5pt]
&\indent+\frac{[a_1,\cdots,a_r]_\infty}{[b_1,\cdots,b_s]_\infty}
\sum_{n=-\infty\atop n\neq0}^\infty
\frac{(-1)^{(s-r+1)n}q^{(s-r+1)n(n+1)/2-n}}{1-q^n}
\left(\frac{a_1\cdots a_r}{b_1\cdots b_s}\right)^n
\nonumber\\[5pt]
&=\frac{[a_1/b_1,\cdots,a_r/b_1]_\infty}{[b_1,b_2/b_1,\cdots,b_s/b_1]_\infty}\\
&\indent\indent\indent\indent\times\sum_{n=-\infty}^\infty\frac{(-1)^{(s-r+1)n}q^{(s-r+1)n(n+1)/2-n}}{1-b_1q^n}
\left(\frac{a_1\cdots a_rb_1^{s-r}}{b_2\cdots b_s}\right)^n
\nonumber\\[5pt]
&\indent\indent\indent\indent\indent\indent+\mathrm{idem}(b_1;b_2,\ldots,b_s).
\end{align}
For $r=s+1$, this is true provided that $|q^2|<|\frac{a_1\cdots a_r}{b_1\cdots b_s}|<|q|$.
\end{thm}

\pf
The proof is similar to that of Theorem \ref{main-result}.
Denote the term on the left-hand side of \eqref{Chan-2.1} by $L'$, and those on the right-hand side by $R'_0,\ldots,R'_s$ respectively, which is
\[L'=R'_0+R'_1+\cdots+R'_s.\]
Likewise by taking $b_0\rightarrow q$ directly in terms other than $R'_0$, we get the right-hand side in \eqref{R'}.

The difference arises in $R'_0$. The terms with $n\neq-1$ are no longer vanishing while taking $b_0\rightarrow q$, which results in an extra generalized Lambert series.
In this case we have
\begin{align*}
&\lim_{b_0\rightarrow q}R'_0
=\lim_{b_0\rightarrow q}\frac{1}{b_0-q}\frac{[a_1/b_0,\cdots,a_r/b_0]_\infty}
{[b_1/b_0,\cdots,b_s/b_0]_\infty}\cdot\frac{(-b_0)^{r-s}b_1\cdots b_s}{a_1\cdots a_r}\cdot q
\nonumber\\[5pt]
&+\frac{[a_1/q,\cdots,a_r/q]_\infty}
{[b_1/q,\cdots,b_s/q]_\infty}
\sum_{n=-\infty\atop n\neq-1}^\infty\frac{(-1)^{(s-r+1)n}q^{(s-r+1)n(n+1)/2+(s-r)n}}{1-q^{n+1}}
\left(\frac{a_1\cdots a_r}{b_1\cdots b_s}\right)^n.
\end{align*}

Thus, denoting the first term by $R''_0$, it suffices to show
\begin{align}
\lim_{b_0\rightarrow q}(L'-R''_0)&=\frac{[a_1,\cdots,a_r]_\infty}
{[b_1,\cdots,b_s]_\infty}\left(1-\mathcal{S}(\bfa)+\mathcal{S}(\bfb)\right).
\end{align}
This can be proved following similar procedures in proving \eqref{L-R_0}.
\qed

\section{An Algorithm for $\mathcal{S}$-series}

The generalized Lambert series $\mathcal{S}$ defined in (\ref{F(a,q)}) appears as encumbrance in our expansions for infinite products.
In this section, we show that $\mathcal{S}(\pm q^m;q^n)$ with $m,n$ integers can be expanded as sums of infinite products.
Therefore, our main results Theorem \ref{main-result} and \ref{single-pole} establish a bridge between infinite products and generalized Lambert series.
We first give some trivial properties concerning $\mathcal{S}$. The following lemma shows that, for special $\bfa$, the function $\mathcal{S}(\bfa)$  degenerates to concise forms.

\begin{lem}\label{Lem-F}
The function $\mathcal{S}$ has the following properties:
\begin{enumerate}
  \item $\mathcal{S}(-1)=-\frac{1}{2}$, $\mathcal{S}(-q)=\frac{1}{2}${\rm;}
  \item $\mathcal{S}(aq)=\mathcal{S}(a)+1${\rm;}
  \item $\mathcal{S}(q/a)=-\mathcal{S}(a)${\rm;}
  \item Let $\bfa=(a_1,\ldots,a_r)$. If $(q/a_1,\ldots,q/a_r)$ is a permutation of $\bfa$, we have $\mathcal{S}(\bfa)=0${\rm;}
  \item $\mathcal{S}(q^s;q^{-t})=\mathcal{S}(q^{s+t};q^t)$.
\end{enumerate}
\end{lem}

In view of (2) and (5), it suffices to consider $\mathcal{S}(\pm q^m;q^n)$ with $m,n$ positive integers.
The proof is trivial, though one should be scrupulous in considering the order of summations in \eqref{F(a,q)}.

\pf
(1) According to the definition of $\mathcal{S}(\bfa)$, we obtain
\begin{align*}
\mathcal{S}(-1)&=\lim_{m\rightarrow\infty}\sum_{n=0}^m\left(\frac{1}{1+q^n}-\frac{1}{1+q^{n+1}}\right)\\
&=\lim_{m\rightarrow\infty}\left(\frac{1}{2}-\frac{1}{1+q^{m+1}}\right)
=-\frac{1}{2}.
\end{align*}
Consequently by (2), we have
$$\mathcal{S}(-q)=\mathcal{S}(-1)+1=\frac{1}{2}.$$

(2) Similarly, we have
\begin{align*}
&\mathcal{S}(aq)-\mathcal{S}(a)
\\[3pt]
&=\lim_{m\rightarrow\infty}\sum_{n=0}^m\left(\frac{1}{1-aq^{n+1}}-\frac{1}{1-q^n/a}\right)
-\lim_{m\rightarrow\infty}\sum_{n=0}^m\left(\frac{1}{1-aq^{n}}-\frac{1}{1-q^{n+1}/a}\right)
\\[3pt]
&=\lim_{m\rightarrow\infty}\left(\frac{1}{1-q^{m+1}/a}-\frac{1}{1-a}+\frac{1}{1-aq^{m+1}}-\frac{1}{1-1/a}\right)
=1.
\end{align*}

(3) By definition, we have
\begin{align*}
\mathcal{S}(q/a)&=\sum_{n=0}^\infty\left(\frac{1}{1-q^{n+1}/a}-\frac{1}{1-aq^n}\right)
\\[3pt]
&=-\sum_{n=0}^\infty\left(\frac{1}{1-aq^n}-\frac{1}{1-q^{n+1}/a}\right)
=-\mathcal{S}(a).
\end{align*}

(4) This follows directly by (3).

(5) We have
\begin{align*}
\mathcal{S}(q^s;q^{-t})&=\sum_{n=0}^\infty\left(\frac{1}{1-q^sq^{-tn}}-\frac{1}{1-q^{-s}q^{-tn-t}}\right)
\\[3pt]
&=\sum_{n=0}^\infty\frac{q^{s-tn}-q^{-s-tn-t}}{(1-q^{s-tn})(1-q^{-s-tn-t})}.
\end{align*}
Multiply both the denominator and numerator of each term by $q^{s+tn+t}q^{-s+tn}$, we derive
\begin{align*}
\mathcal{S}(q^s;q^{-t})&=\sum_{n=0}^\infty\left(\frac{1}{1-q^{s+tn+t}}-\frac{1}{1-q^{-s+tn}}\right)
=\mathcal{S}(q^{s+t};q^{t}).
\end{align*}
\qed

The following lemma is due to Andrews, Lewis and Liu \cite{Andrews-Lewis-Liu-2001}. Chan \cite{Chan-2005} provided another proof using Lemma \ref{Chan-Thm-2.2}.
\begin{lem}\label{Chan-Coro-3.2} For $|q|<1$, we have
\begin{align}
\frac{[ab,bc,ca]_\infty(q)_\infty^2}{[a,b,c,abc]_\infty}
&=1+\sum_{n=0}^\infty\frac{aq^n}{1-aq^n}-\sum_{n=1}^\infty\frac{q^n/a}{1-q^n/a}
+\sum_{n=0}^\infty\frac{bq^n}{1-bq^n}
\nonumber\\[5pt]
&\quad-\sum_{n=1}^\infty\frac{q^n/b}{1-q^n/b}
+\sum_{n=0}^\infty\frac{cq^n}{1-cq^n}-\sum_{n=1}^\infty\frac{q^n/c}{1-q^n/c}
\nonumber\\[5pt]
&\quad-\sum_{n=0}^\infty\frac{abcq^n}{1-abcq^n}
+\sum_{n=1}^\infty\frac{q^n/abc}{1-q^n/abc}.
\end{align}
\end{lem}
Lemma \ref{Chan-Coro-3.2} associates the function $\mathcal{S}$ with theta functions.
We denote the infinite products on the left-hand side of Lemma \ref{Chan-Coro-3.2} by $\mathcal{P}(a,b,c)$, which is
\begin{equation*}
\mathcal{P}(a,b,c)=\mathcal{P}(a,b,c;q)=\frac{[ab,bc,ca]_\infty(q)_\infty^2}{[a,b,c,abc]_\infty}.
\end{equation*}
For the sake of brevity, we denote $\mathcal{P}(a,a,a)$ by $\mathcal{P}(a)$.
Then, Lemma \ref{Chan-Coro-3.2} shows that
\begin{align}\label{Chan-Coro-3.2-huajian}
\mathcal{P}(a,b,c)=1+\mathcal{S}(a)
+\mathcal{S}(b)+\mathcal{S}(c)-\mathcal{S}(abc).
\end{align}

We are now equipped to propose an algorithm for $\mathcal{S}(\pm q^m;q^n)$ with arbitrary positive integers $m$ and $n$.
First in \eqref{Chan-Coro-3.2-huajian}, by replacing $q$ by $q^n$ and setting $a=\pm q^m$, $b=\pm q^m$ and $c=-q^{n-2m}$, we have
\begin{align}\label{G(q^m,q^m,-q^{n-2m}}
\nonumber&\mathcal{P}(\pm q^m,\pm q^m,-q^{n-2m};q^n)\\
\nonumber&\indent\indent=1+2\mathcal{S}(\pm q^m;q^n)
+\mathcal{S}(-q^{n-2m};q^n)
-\mathcal{S}(-q^n;q^n)
\\
&\indent\indent=\frac{1}{2}+2\mathcal{S}(\pm q^m;q^n)-\mathcal{S}(-q^{2m};q^n).
\end{align}
Therefore, in order to obtain expansions for $\mathcal{S}(\pm q^m;q^n)$ in terms of $\mathcal{P}$-functions, we need to calculate $\mathcal{S}(-q^{2m};q^n)$.
Our strategy is to implement a recursive procedure using (\ref{Chan-Coro-3.2-huajian}).

Suppose that $n=3^s\cdot n'$ with $(3,n')=1$. We denote by $k$ the order of $3$ in the cyclic group $\mathbb{Z}_{n'}$, which is
\begin{equation}\label{k}
k=k(n')=\mathrm{ord}_{\mathbb{Z}_{n'}}(3).
\end{equation}
Thus, we have
\begin{equation*}
3^k\equiv 1 \mmod{n'}
\end{equation*}
and accordingly
\begin{equation}\label{3^{s+k}}
3^{s+k}\equiv 3^{s} \mmod{n}.
\end{equation}
Then, by setting all $a,b,c$ with $-q^{3^{j-1}\cdot 2m}$ where $j=1,\ldots,s+k$ successively, we obtain a chain of identities as following:
\begin{align*}
j=1:&&\mathcal{P}(-q^{2m};q^n)&=1+3\mathcal{S}(-q^{2m};q^n)-\mathcal{S}(-q^{3\cdot 2m};q^n);\\
\vdots\indent~~&&\vdots\indent~~&\indent\indent\indent\indent\vdots\\
j=s:&&\mathcal{P}(-q^{3^{s-1}\cdot 2m};q^n)&=1+3\mathcal{S}(-q^{3^{s-1}\cdot 2m};q^n)-\mathcal{S}(-q^{3^{s}\cdot 2m};q^n);\\
j=s+1:&&\mathcal{P}(-q^{3^{s}\cdot 2m};q^n)&=1+3\mathcal{S}(-q^{3^{s}\cdot 2m};q^n)-\mathcal{S}(-q^{3^{s+1}\cdot 2m};q^n);\\
\vdots\indent~~&&\vdots\indent~~&\indent\indent\indent\indent\vdots\\
j=s+k:&&\mathcal{P}(-q^{3^{s+k-1}\cdot 2m};q^n)&=1+3\mathcal{S}(-q^{3^{s+k-1}\cdot 2m};q^n)-\mathcal{S}(-q^{3^{s+k}\cdot 2m};q^n).
\end{align*}
In view of \eqref{3^{s+k}} and Lemma \ref{Lem-F}(2), we are now able to solve $\mathcal{S}(-q^{2m};q^n)$.
Concretely, for $j=1,\ldots,s$, we multiply the identities by $3^{s-j}$ respectively. Then, their weighted summation turns to
\begin{equation}\label{s-identities}
\sum_{j=1}^{s}3^{s-j}\mathcal{P}(-q^{3^{j-1}\cdot 2m};q^n)
=\frac{3^s-1}{2}+3^s\mathcal{S}(-q^{2m};q^n)-\mathcal{S}(-q^{3^{s}\cdot 2m};q^n).
\end{equation}
Again for $j=s+1,\ldots,s+k$, we multiply the identities by $3^{s+k-j}$ respectively. Then, their weighted summation turns to
\begin{equation}\label{k-identities}
\sum_{j=s+1}^{s+k}3^{s+k-j}\mathcal{P}(-q^{3^{j-1}\cdot 2m};q^n)
=\frac{3^k-1}{2}+3^k\mathcal{S}(-q^{3^{s}\cdot 2m};q^n)-\mathcal{S}(-q^{3^{s+k}\cdot 2m};q^n).
\end{equation}
Considering \eqref{3^{s+k}}, we have
\begin{equation}\label{3^(s+k)*2m}
\mathcal{S}(-q^{3^{s+k}\cdot 2m};q^n)=\mathcal{S}(-q^{3^{s}\cdot 2m};q^n)+\frac{3^s(3^k-1)\cdot 2m}{n}.
\end{equation}
Combining \eqref{s-identities}, \eqref{k-identities} and \eqref{3^(s+k)*2m}, we are able to obtain $\mathcal{S}(-q^{2m};q^n)$, and consequently $\mathcal{S}(\pm q^m;q^n)$ by \eqref{G(q^m,q^m,-q^{n-2m}}. We summarize the algorithm as the following theorem.

\begin{thm}\label{F-main}
Suppose that $m$ and $n$ are positive integers with $n=3^s\cdot n'$ and $(3,n')=1$. Denote by $k$ the order of $3$ in the cyclic group $\mathbb{Z}_{n'}$. Then, we have
\begin{align*}
2\mathcal{S}(\pm q^m;q^n)
&+\frac{n-2m}{n}\\[3pt]
&=\sum_{j=1}^{s+k}\frac{3^{k-j}}{3^k-1}\mathcal{P}(-q^{3^{j-1}\cdot 2m};q^n)-\sum_{j=1}^{s}\frac{3^{-j}}{3^k-1}\mathcal{P}(-q^{3^{j-1}\cdot 2m};q^n)\\[5pt]
&\indent\indent\indent\indent\indent+\mathcal{P}(\pm q^m,\pm q^m,-q^{n-2m};q^n).
\end{align*}
\end{thm}

The length of the chain may be reduced for special $m$ and $n$.
We consider the first $l$ identities in the chain. Their summation with weights $3^{l-j}$ gives
\begin{equation}\label{l-identities}
\sum_{j=1}^{l}3^{l-j}\mathcal{P}(-q^{3^{j-1}\cdot 2m};q^n)
=\frac{3^l-1}{2}+3^l\mathcal{S}(-q^{2m};q^n)-\mathcal{S}(-q^{3^{l}\cdot 2m};q^n).
\end{equation}
Lemma \ref{Lem-F} provides values of $\mathcal{S}$ at special points, which would help to shorten the chain of identities.
Suppose that
\begin{align*}
&n=3^{s_1}\cdot 2^{t_1}\cdot n' \indent\mathrm{with} \indent (3,n')=1 ~\mathrm{and}~ (2,n')=1,\\
&m=3^{s_2}\cdot 2^{t_2}\cdot m' \indent\mathrm{with} \indent (3,m')=1 ~\mathrm{and}~ (2,m')=1.
\end{align*}
We consider two special cases.

\emph{Case I:} $n'\mid m'$ and $t_1\le t_2+1$.

We take $l$ by setting
\begin{equation}\label{l}
l=\begin{cases} 0,&\text{when $s_2\ge s_1$,}\\
s_1-s_2,&\text{when $s_2< s_1$.}\end{cases}
\end{equation}
In this case, $l$ is the least nonnegative integer such that
$$3^l\cdot 2m\equiv 0\mmod n.$$
By Lemma \ref{Lem-F}, we have
\begin{align}\label{3^l*2m-1}
\nonumber\mathcal{S}(-q^{3^{l}\cdot 2m};q^n)
&=\frac{3^l\cdot 2m}{n}+\mathcal{S}(-1;q^n)\\
&=\frac{3^l\cdot 2m}{n}-\frac{1}{2}.
\end{align}
Combining \eqref{l-identities} and \eqref{3^l*2m-1}, we are able to obtain $\mathcal{S}(-q^{2m};q^n)$, and consequently $\mathcal{S}(q^m;q^n)$ by \eqref{G(q^m,q^m,-q^{n-2m}}.

\emph{Case II:} $n'\mid m'$ and $t_1=t_2+2$.

We take $l$ as in \eqref{l}.
Now $l$ is the least nonnegative integer such that
$$3^l\cdot 2m\equiv n/2\mmod n.$$
The discussion is similar to that of Case I. A tiny difference lies in \eqref{3^l*2m-1}, where we now have
\begin{align}\label{3^l*2m-2}
\nonumber\mathcal{S}(-q^{3^{l}\cdot 2m};q^n)
&=\frac{3^l\cdot 2m-n/2}{n}+\mathcal{S}(q^\frac{n}{2};q^n)\\
&=\frac{3^l\cdot 2m}{n}-\frac{1}{2}.
\end{align}

We summarize these two cases as the following corollary.
\begin{coro}\label{cor-F_1}
Let $m$ and $n$ be positive integers. Suppose that there exists a least nonnegative integer $l$ such that $3^l\cdot 4m\equiv 0\mmod n$. Then, we have
\begin{align*}
2\mathcal{S}(\pm q^m;q^n)&+\frac{n-2m}{n}\\
&=\sum_{j=1}^{l}3^{-j}\mathcal{P}(-q^{3^{j-1}\cdot 2m};q^n)
+\mathcal{P}(\pm q^m,\pm q^m,-q^{n-2m};q^n).
\end{align*}
\end{coro}

For example, when $n=3$, we have $l=1$ in Corollary \ref{cor-F_1}. We give the explicit expansion for $\mathcal{S}(\pm q;q^3)$ in terms of infinite products.
\begin{coro}
We have
\begin{align}\label{Chan-Coro-gene}
\nonumber1+6\mathcal{S}(q;q^3)
&=3\mathcal{P}(q,q,-q;q^3)
-\mathcal{P}(-q,-q,-q;q^3)
\\[5pt]
&=\frac{3J_{3,6}^2J_{6}^3}{2J_{1,6}^2J_2}
-\frac{J_{1,6}^6J_{2}^3J_{3,6}^2}{2J_{6}^9}
\end{align}
and
\begin{align}\label{Coro-gene}
1+6\mathcal{S}(-q;q^3)
=2\mathcal{P}(-q,-q,-q;q^3)
=\frac{J_{1,6}^6J_{2}^3J_{3,6}^2}{J_{6}^9}.
\end{align}
\end{coro}

\section{Examples for 3-dissections}
Combined with the algorithm of $\mathcal{S}$-series, Theorem \ref{main-result} and \ref{single-pole} constructed a bridge between sums of generalized Lambert series and those of theta functions.
In this section, We show some examples concerning 3-dissections of generalized Lambert series.
These formulas make comparison between Chan's identities and ours in this article.
In \S5, they are useful to discuss properties of ranks of overpartitions modulo 6.

First, by both replacing $q$ by $q^3$ and taking $b_0=q$, $b_1=-q$ in Lemma \ref{Chan-Thm-2.1} and Theorem \ref{main-result} respectively, we have the following corollary.
This shows that, the double poles make it more complicated for the correspondence between generalized Lambert series and infinite products.
\begin{coro}\label{cor4.1}
We have
\begin{align}
&\sum_{n=-\infty}^\infty\frac{(-1)^n q^{3n^2+3n}}{1+q^{3n+1}}
+\sum_{n=-\infty}^\infty\frac{(-1)^n q^{3n^2+3n}}{1-q^{3n+1}}
=\frac{2J^3_{6}}{J_{2}},\label{simpli1}\\
&\sum_{n=-\infty}^\infty\frac{(-1)^n q^{3n^2+3n}}{(1+q^{3n+1})^2}
+\sum_{n=-\infty}^\infty\frac{(-1)^n q^{3n^2+3n}}{(1-q^{3n+1})^2}
=\frac{4J^3_{6}}{3J_{2}}+\frac{J_{3,6}^2J_{6}^6}{2J_{1,6}^2J_2^2}+\frac{J_{1,6}^6J_{2}^2J_{3,6}^2}{6J_{6}^6}.\label{simpli2}
\end{align}
\end{coro}

The following corollary makes comparison between Lemma \ref{Chan-Thm-2.1} and Theorem \ref{single-pole}.
\begin{coro}\label{cor3.6}
We have
\begin{align}
\nonumber\sum_{n=-\infty}^\infty\frac{(-1)^nq^{9n^2+3n}}{1+q^{9n}}
&+\sum_{n=-\infty}^\infty\frac{(-1)^nq^{9n^2+15n+6}}{1+q^{9n+6}}\\
&=\frac{2J_{3,18}}{J_{9,18}}
\sum_{n=-\infty}^\infty\frac{(-1)^nq^{9n^2+9n+3}}{1+q^{9n+3}}
+\frac{ J_{3,18}^6J_{6}^3J_{9,18}^2}{2J_{18}^9},
\label{Chan-Thm-2.1-coro}\\
\nonumber\sum_{n=-\infty\atop n\neq0}^{\infty}\frac{(-1)^nq^{9n^2+3n}}{1-q^{9n}}
&+\sum_{n=-\infty}^{\infty}\frac{(-1)^nq^{9n^2+15n+6}}{1-q^{9n+6}}\\
&=\frac{2J_{3,18}}{J_{9,18}}\sum_{n=-\infty}^\infty
\frac{(-1)^nq^{9n^2+9n+3}}{1-q^{9n+3}}
+\frac{J_{3,18}^6J_{6}^3J_{9,18}^2}{6J_{18}^9}
-\frac{1}{6}.\label{main-result-coro}
\end{align}
\end{coro}

\pf
For \eqref{Chan-Thm-2.1-coro}, we set $r=1$ and $s=3$ in Lemma \ref{Chan-Thm-2.1}. By replacing $q$ by $q^9$ and taking $a_1=q^3$, $b_1=-1$, $b_2=-q^3$, and $b_3=-q^6$, we obtain
\begin{align*}
\sum_{n=-\infty}^{\infty}\frac{(-1)^nq^{9n^2+3n}}{1+q^{9n}}
&+\sum_{n=-\infty}^{\infty}\frac{(-1)^nq^{9n^2+15n+6}}{1+q^{9n+6}}
\nonumber\\[5pt]
&=\frac{[-1;q^9]_\infty}{[-q^3;q^9]_\infty}\sum_{n=-\infty}^\infty
\frac{(-1)^nq^{9n^2+9n+3}}{1+q^{9n+3}}
+\mathcal{P}(-q^3;q^9).
\end{align*}

For \eqref{main-result-coro}, we set $r=1$ and $s=2$ in Theorem \ref{single-pole}. By replacing $q$ by $q^9$ and taking $a_1=-q^{12}$, $b_1=q^3$ and $b_2=q^6$, we obtain
\begin{align*}
\sum_{n=-\infty\atop n\neq0}^{\infty}\frac{(-1)^nq^{9n^2+3n}}{1-q^{9n}}
&+\sum_{n=-\infty}^{\infty}\frac{(-1)^nq^{9n^2+15n+6}}{1-q^{9n+6}}
\nonumber\\[5pt]
&=\frac{[-1;q^9]_\infty}{[-q^3;q^9]_\infty}\sum_{n=-\infty}^\infty
\frac{(-1)^nq^{9n^2+9n+3}}{1-q^{9n+3}}
+\mathcal{S}(-q^3;q^9).
\end{align*}
Then \eqref{main-result-coro} follows by \eqref{Coro-gene}.
\qed

Consider 3-dissections of generalized Lambert series according to the summation index $n$ modulo 3:
\begin{align*}
\sum_{n=-\infty}^\infty\frac{(-1)^nq^{n^2+n}}{1+q^{3n}}
&=\sum_{n=-\infty}^\infty\frac{(-1)^nq^{9n^2+3n}}{1+q^{9n}}\\
&\indent\indent
-\sum_{n=-\infty}^\infty\frac{(-1)^nq^{9n^2+9n+2}}{1+q^{9n+3}}
+\sum_{n=-\infty}^\infty\frac{(-1)^nq^{9n^2+15n+6}}{1+q^{9n+6}},\\
\sum_{n=-\infty\atop n\neq0}^\infty\frac{(-1)^nq^{n^2+n}}{1-q^{3n}}
&=\sum_{n=-\infty\atop n\neq0}^\infty\frac{(-1)^nq^{9n^2+3n}}{1-q^{9n}}\\
&\indent\indent
-\sum_{n=-\infty}^\infty\frac{(-1)^nq^{9n^2+9n+2}}{1-q^{9n+3}}
+\sum_{n=-\infty}^\infty\frac{(-1)^nq^{9n^2+15n+6}}{1-q^{9n+6}}.
\end{align*}
Using Corollary \ref{cor3.6}, we transform these 3-dissections into forms containing one single generalized Lambert series.

\begin{coro}\label{cor3.7}
We have
\begin{align*}
&\sum_{n=-\infty}^\infty\frac{(-1)^nq^{n^2+n}}{1+q^{3n}}\\
&\indent\indent=\left(2q\frac{J_{3,18}}{J_{9,18}}-1\right)
\sum_{n=-\infty}^\infty\frac{(-1)^nq^{9n^2+9n+2}}{1+q^{9n+3}}
+\frac{J_{3,18}^6J_{6}^3J_{9,18}^2}{2J_{18}^9},\\
&\sum_{n=-\infty\atop n\neq0}^\infty\frac{(-1)^nq^{n^2+n}}{1-q^{3n}}\\
&\indent\indent=\left(2q\frac{J_{3,18}}{J_{9,18}}-1\right)
\sum_{n=-\infty}^\infty
\frac{(-1)^nq^{9n^2+9n+2}}{1-q^{9n+3}}
+\frac{J_{3,18}^6J_{6}^3J_{9,18}^2}{6J_{18}^9}
-\frac{1}{6}.
\end{align*}
\end{coro}

As we mentioned, Theorem \ref{main-result} decouples parameters in $\bfa$ from poles in generalized Lambert series.
This helps in constructing identities with variant orders in the numerators.
Consider the following 3-dissections of generalized Lambert series:
\begin{align*}
\sum_{n=-\infty}^\infty\frac{(-1)^nq^{n^2+n}}{(1+q^{3n})^2}
&=\sum_{n=-\infty}^\infty\frac{(-1)^nq^{9n^2+3n}}{(1+q^{9n})^2}\\
&\indent\indent
-\sum_{n=-\infty}^\infty\frac{(-1)^nq^{9n^2+9n+2}}{(1+q^{9n+3})^2}
+\sum_{n=-\infty}^\infty\frac{(-1)^nq^{9n^2+15n+6}}{(1+q^{9n+6})^2},\\
\sum_{n=-\infty}^\infty\frac{(-1)^nq^{n^2+2n}}{(1+q^{3n})^2}
&=\sum_{n=-\infty}^\infty\frac{(-1)^nq^{9n^2+6n}}{(1+q^{9n})^2}\\
&\indent\indent
-\sum_{n=-\infty}^\infty\frac{(-1)^nq^{9n^2+12n+3}}{(1+q^{9n+3})^2}
+\sum_{n=-\infty}^\infty\frac{(-1)^nq^{9n^2+18n+8}}{(1+q^{9n+6})^2},\\
\sum_{n=-\infty}^\infty\frac{(-1)^nq^{n^2+3n}}{(1+q^{3n})^2}
&=\sum_{n=-\infty}^\infty\frac{(-1)^nq^{9n^2+9n}}{(1+q^{9n})^2}\\
&\indent\indent
-\sum_{n=-\infty}^\infty\frac{(-1)^nq^{9n^2+15n+4}}{(1+q^{9n+3})^2}
+\sum_{n=-\infty}^\infty\frac{(-1)^nq^{9n^2+21n+10}}{(1+q^{9n+6})^2}.
\end{align*}
Similar to Corollary \ref{cor3.7}, we aim to transform these 3-dissections into forms containing one single generalized Lambert series.
\begin{coro}\label{cor3.8}
We have
\begin{align*}
\sum_{n=-\infty}^\infty\frac{(-1)^nq^{n^2+n}}{(1+q^{3n})^2}
&=\frac{J_{3,18}^6J_6^3J_{9,18}^2}{2J_{18}^9}
\left(\frac{2}{3}-\frac{J_{9,18}^2J_{18}^3}{4J_{3,18}^2J_6}
+\frac{J_{3,18}^6J_{6}^3J_{9,18}^2}{12J_{18}^9}\right)\\
&\indent\indent\indent\indent-\left(1-\frac{2qJ_{3,18}}{J_{9,18}}\right)
\sum_{n=-\infty}^\infty\frac{(-1)^nq^{9n^2+9n+2}}{(1+q^{9n+3})^2},\\
\sum_{n=-\infty}^\infty\frac{(-1)^nq^{n^2+2n}}{(1+q^{3n})^2}
&=\frac{J_{3,18}^6J_6^3J_{9,18}^2}{2J_{18}^9}
\left(\frac{1}{3}+\frac{J_{9,18}^2J_{18}^3}{4J_{3,18}^2J_6}
-\frac{J_{3,18}^6J_{6}^3J_{9,18}^2}{12J_{18}^9}\right)\\
&\indent\indent\indent\indent-\left(1-\frac{2qJ_{3,18}}{J_{9,18}}\right)
\sum_{n=-\infty}^\infty\frac{(-1)^nq^{9n^2+18n+5}}{(1+q^{9n+3})^2},\\
\sum_{n=-\infty}^\infty\frac{(-1)^nq^{n^2+3n}}{(1+q^{3n})^2}
&=\frac{qJ_{3,18}^5J_{6}^2J_{9,18}^3}{2J_{18}^6}
+\left(1-\frac{2qJ_{3,18}}{J_{9,18}}\right)
\sum_{n=-\infty}^\infty\frac{(-1)^nq^{9n^2+9n}}{(1+q^{9n})^2}.
\end{align*}
\end{coro}

\pf
For the third identity, we replace $q$ by $q^9$ and set $r=1$, $s=3$ in Lemma \ref{Chan-Thm-3.2}.
Then by taking $a_0=1$, $b_1=-1$, $b_2=-q^3$, $b_3=-q^6$, we obtain
\begin{align}
\nonumber-\sum_{n=-\infty}^\infty\frac{(-1)^nq^{9n^2+15n+3}}{(1+q^{9n+3})^2}
&+\sum_{n=-\infty}^\infty\frac{(-1)^nq^{9n^2+21n+9}}{(1+q^{9n+6})^2}\\
&=\frac{J_{6}^2J_{3,18}^5J_{9,18}^3}{2J_{18}^6}
-\frac{2J_{3,18}}{J_{9,18}}
\sum_{n=-\infty}^\infty\frac{(-1)^nq^{9n^2+9n}}{(1+q^{9n})^2}.\label{q^{n^2+3n}}
\end{align}
This proves the third identity in the corollary.

When concerning the first and second identities, Lemma \ref{Chan-Thm-3.2} fails to give a proper relationship similar to \eqref{q^{n^2+3n}}.
Poles are twisted with the parameter $a_0$, which limits the orders of $q$ in numerators.
Instead we replace $q$ by $q^9$ and set $r=1$, $s=3$ in Theorem \ref{main-result}.
Then by taking $a_1=q^3$, $b_1=-1$, $b_2=-q^3$, $b_3=-q^6$, we obtain
\begin{align*}
&\sum_{n=-\infty}^\infty\frac{(-1)^nq^{9n^2+3n}}{(1+q^{9n})^2}
+\sum_{n=-\infty}^\infty\frac{(-1)^nq^{9n^2+15n+6}}{(1+q^{9n+6})^2}\\
&=\frac{J_6^3J_{3,18}^6J_{9,18}^2}{2J_{18}^9}
\left(\frac{1}{2}-\mathcal{S}(q^3;q^9)\right)
+\frac{2J_{3,18}}{J_{9,18}}
\sum_{n=-\infty}^\infty\frac{(-1)^nq^{9n^2+9n+3}}{(1+q^{9n+3})^2}.
\end{align*}
Then the second identity follows by \eqref{Chan-Coro-gene}.

Similarly, by taking $a_1=q^6$, $b_1=-1$, $b_2=-q^3$, $b_3=-q^6$, we obtain
\begin{align*}
&\sum_{n=-\infty}^\infty\frac{(-1)^nq^{9n^2+6n}}{(1+q^{9n})^2}
-\sum_{n=-\infty}^\infty\frac{(-1)^nq^{9n^2+12n+3}}{(1+q^{9n+3})^2}\\
&=\frac{J_6^3J_{3,18}^6J_{9,18}^2}{2J_{18}^9}
\left(\frac{1}{2}-\mathcal{S}(q^6;q^9)\right)
-\frac{2J_{3,18}}{J_{9,18}}
\sum_{n=-\infty}^\infty\frac{(-1)^nq^{9n^2+18n+9}}{(1+q^{9n+6})^2}.
\end{align*}
Noting that
\begin{align*}
\sum_{n=-\infty}^\infty\frac{(-1)^nq^{9n^2+18n+8}}{(1+q^{9n+6})^2}
=-\sum_{n=-\infty}^\infty\frac{(-1)^nq^{9n^2+18n+5}}{(1+q^{9n+3})^2}.
\end{align*}
Then the second identity follows by \eqref{Chan-Coro-gene} and
\begin{equation*}
\mathcal{S}(q^6;q^9)=-\mathcal{S}(q^3;q^9).
\end{equation*}
\qed

\section{Ranks of Overpartitions modulo $6$}
In this section, we study 3-dissection properties of ranks of overpartitions modulo 6.
Noting that
\begin{equation*}
\overline{N}(s,\ell,n)=\overline{N}(\ell-s,\ell,n),
\end{equation*}
it suffices to consider four residue classes when $n=6$.

Replacing $z$ by $\xi_6=e^{\frac{\pi i}{3}}$, the root of unity modulo $6$, left-hand side of \eqref{GF} reduces to
\begin{align*}
\overline{R}(\xi_6;q)
&=\sum_{n=0}^\infty\sum_{m=-\infty}^{\infty}\overline{N}(m,n)\xi_6^mq^n
\nonumber\\[5pt]
&=\sum_{n=0}^\infty\sum_{t=0}^5\sum_{m=-\infty}^{\infty}\overline{N}(6m+t,n)\xi_6^tq^n
\nonumber\\[5pt]
&=\sum_{n=0}^\infty(\overline{N}(0,6,n)+\overline{N}(1,6,n)
-\overline{N}(2,6,n)-\overline{N}(3,6,n))q^n.
\end{align*}
On the other hand, in light of the fact that $\xi_6+\xi_6^{-1}=1$, we get
\begin{align*}
\overline{R}(\xi_6;q)
&=\frac{(-q)_\infty}{(q)_\infty}\left\{1+2\sum_{n=1}^\infty
\frac{(2-\xi_6-\xi_6^{-1})(-1)^nq^{n^2+n}}{1-\xi_6q^n-\xi_6^{-1}q^n+q^{2n}}\right\}\\
&=\frac{(-q)_\infty}{(q)_\infty}\left\{1+2\sum_{n=1}^\infty
\frac{(-1)^nq^{n^2+n}}{1-q^n+q^{2n}}\right\}\\
&=\frac{2(-q)_\infty}{(q)_\infty}\sum_{n=-\infty}^\infty
\frac{(-1)^nq^{n^2+n}}{1+q^{3n}}.
\end{align*}
Thus, we get
\begin{align}\label{xi_6}
\nonumber\overline{R}(\xi_6;q)
&=\sum_{n=0}^\infty(\overline{N}(0,6,n)+\overline{N}(1,6,n)
-\overline{N}(2,6,n)-\overline{N}(3,6,n))q^n\\
&=\frac{2(-q)_\infty}{(q)_\infty}\sum_{n=-\infty}^\infty
\frac{(-1)^nq^{n^2+n}}{1+q^{3n}}.
\end{align}
Similarly, if we replace $z$ by $\xi_6^2$, $\xi_6^3$ and $1$ in the left-hand side of \eqref{GF} respectively, we obtain
\begin{align}
\nonumber\overline{R}(\xi_6^2;q)
&=\sum_{n=0}^\infty(\overline{N}(0,6,n)-\overline{N}(1,6,n)
-\overline{N}(2,6,n)+\overline{N}(3,6,n))q^n\\
&=\frac{6(-q)_\infty}{(q)_\infty}
\sum_{n=-\infty\atop n\neq0}^\infty\frac{(-1)^nq^{n^2+n}}{1-q^{3n}}+\frac{(-q)_\infty}{(q)_\infty},\label{xi_6^2}\\
\nonumber\overline{R}(\xi_6^3;q)
&=\sum_{n=0}^\infty(\overline{N}(0,6,n)-2\overline{N}(1,6,n)
+2\overline{N}(2,6,n)-\overline{N}(3,6,n))q^n\\
&=\frac{4(-q)_\infty}{(q)_\infty}\sum_{n=-\infty}^\infty
\frac{(-1)^nq^{n^2+n}}{(1+q^{n})^2},\label{xi_6^3}\\
\nonumber\overline{R}(1;q)
&=\sum_{n=0}^\infty(\overline{N}(0,6,n)+2\overline{N}(1,6,n)
+2\overline{N}(2,6,n)+\overline{N}(3,6,n))q^n\\
&=\frac{(-q)_\infty}{(q)_\infty}.\label{1}
\end{align}
Now, we have a linear equation system concerning all residues of ranks for overpartitions modulo 6.
The rank of its coefficient matrix is full, so we are able to solve $\overline{N}(i,6,n)$ for $i=0,1,2,3$ in terms of $\overline{R}(z;q)$.
\begin{lem}\label{lem4.1} We have
\begin{align*}
\sum_{n=0}^\infty\overline{N}(0,6,n)q^n
&=\frac{1}{6}\left(\overline{R}(1;q)
+2\overline{R}(\xi_6;q)
+2\overline{R}(\xi_6^2;q)
+\overline{R}(\xi_6^3;q)\right),\\[5pt]
\sum_{n=0}^\infty\overline{N}(1,6,n)q^n
&=\frac{1}{6}\left(\overline{R}(1;q)
+~~\overline{R}(\xi_6;q)
-~~\overline{R}(\xi_6^2;q)
-\overline{R}(\xi_6^3;q)\right),\\[5pt]
\sum_{n=0}^\infty\overline{N}(2,6,n)q^n
&=\frac{1}{6}\left(\overline{R}(1;q)
-~~\overline{R}(\xi_6;q)
-~~\overline{R}(\xi_6^2;q)
+\overline{R}(\xi_6^3;q)\right),\\[5pt]
\sum_{n=0}^\infty\overline{N}(3,6,n)q^n
&=\frac{1}{6}\left(\overline{R}(1;q)
-2\overline{R}(\xi_6;q)
+2\overline{R}(\xi_6^2;q)
-\overline{R}(\xi_6^3;q)\right).
\end{align*}
\end{lem}
Therefore, if we can elaborate 3-dissection properties of each $\overline{R}$ function, we can go further to those of overpartitions.
In view of \eqref{xi_6}-\eqref{1}, it is not surprised that the identities we obtained in \S3 will play a key role.
We first give some lemmas.
\begin{lem}\label{Lovejoy-Osburn-2008-mod6}
We have
\begin{align}
\nonumber\frac{(q;q)_\infty}{(-q;q)_\infty}
&=\frac{(q^9;q^9)_\infty}{(-q^9;q^9)_\infty}
-2q(q^3,q^{15},q^{18};q^{18})_\infty\\[5pt]
&=J_{9,18}-2qJ_{3,18}.
\end{align}
\end{lem}
\pf
This is \cite[Theorem 1.2]{Andrews-Hickerson-1991}.
\qed

\begin{lem}\label{V-0}
We have
\begin{equation}\label{V-0-rr}
\frac{J_{18}^3}{J_{3,18}^3}-8q^3\frac{J_{18}^3}{J_{9,18}^3}
=\frac{J_{3,18}^5J_{6}^4}{J_{18}^9}.
\end{equation}
\end{lem}
\pf
This identity is equivalent to
\begin{equation}\label{4.6}
\mathcal{P}(q^3,q^3,-q^3;q^9)+\mathcal{P}(-q^3,-q^3,q^6;q^9)=\mathcal{P}(-q^3,-q^3,-q^3;q^9),
\end{equation}
which can be easily verified by \eqref{Chan-Coro-3.2-huajian}.
In fact, \eqref{4.6} is a special case of the following identity in \cite{Atkin-Swinnerton-Dyer-1954}
\begin{align*}
j(x;q)^2j(yz;q)j(yz^{-1};q)
&=j(y;q)^2j(xz;q)j(xz^{-1};q)\\
&\indent\indent-yz^{-1}j(z;q)^2j(xy;q)j(x y^{-1};q)
\end{align*}
by replacing $q$ by $q^9$, then setting $x=-q^3$, $y=q^3$ and $z=-1$.
\qed

Now, we give the 3-dissections of $\overline{R}(z;q)$ for $z=1,\xi_6,\xi_6^2,\xi_6^3$ successively.

\begin{lem}\label{R(1;q)}
We have
\begin{align*}
&\overline{R}(1;q)
=\frac{J_{18}^{12}}{J_{3,18}^8J_{6}^4J_{9,18}}
+q\frac{2J_{18}^{12}}{J_{3,18}^7J_{6}^4J_{9,18}^2}
+q^2\frac{4J_{18}^{12}}{J_{3,18}^6J_{6}^4J_{9,18}^3}.
\end{align*}
\end{lem}

\pf
Hirschhorn and Sellers \cite{Hirschhorn-Sellers-2005} proved that
\begin{align}\label{Hirschhorn-Sellers-2005}
\frac{(-q)_\infty}{(q)_\infty}
=\frac{J_{18}^{12}}{J_{3,18}^8J_{6}^4J_{9,18}}
+q\frac{2J_{18}^{12}}{J_{3,18}^7J_{6}^4J_{9,18}^2}
+q^2\frac{4J_{18}^{12}}{J_{3,18}^6J_{6}^4J_{9,18}^3}.
\end{align}
Then, the lemma follows by \eqref{1}.
One can also verify \eqref{Hirschhorn-Sellers-2005} easily by Lemma \ref{Lovejoy-Osburn-2008-mod6} and \ref{V-0}.
\qed

\begin{lem}\label{R(xi_6;q)}
We have
\begin{align*}
&\overline{R}(\xi_6;q)
=\frac{J_{18}^3J_{9,18}}{J^2_{3,18}J_{6}}
+q\frac{2J_{18}^3}{J_{3,18}J_{6}}
+q^2\left(\frac{4J^3_{18}}{J_{6}J_{9,18}}
-\frac{2}{J_{9,18}}\sum_{n=-\infty}^\infty\frac{(-1)^n q^{9n^2+9n}}{1+q^{9n+3}}\right),\\
&\overline{R}(\xi_6^2;q)
=\frac{J_{18}^3J_{9,18}}{J_{3,18}^2J_{6}}
+q\frac{2J_{18}^3}{J_{3,18}J_{6}}
+q^2\left(\frac{4J_{18}^3}{J_{6}J_{9,18}}-\frac{6}{J_{9,18}}
\sum_{n=-\infty}^\infty\frac{(-1)^nq^{9n^2+9n}}{1-q^{9n+3}}\right).
\end{align*}
\end{lem}
\pf
Ji, Zhang and Zhao \cite[(2.3)]{Ji-Zhang-Zhao-2017} proved $\overline{R}(\xi_6;q)$ by using  Corollary \ref{cor3.7}.

For $\overline{R}(\xi_6^2;q)$, we substitute the second identity of Corollary \ref{cor3.7} into \eqref{xi_6^2}, and obtain
\begin{align*}
\overline{R}(\xi_6^2;q)
&=\frac{6(-q)_\infty}{(q)_\infty}
\sum_{n=-\infty\atop n\neq0}^\infty\frac{(-1)^nq^{n^2+n}}{1-q^{3n}}+\frac{(-q)_\infty}{(q)_\infty}
\nonumber\\[5pt]
&=-\frac{6}{J_{9,18}}
\sum_{n=-\infty}^\infty\frac{(-1)^nq^{9n^2+9n+2}}{1-q^{9n+3}}
+\frac{(-q)_\infty}{(q)_\infty}
\frac{J_{6}^3J_{3,18}^6J_{9,18}^2}{J_{18}^9}.
\end{align*}
Then, the lemma follows by \eqref{Hirschhorn-Sellers-2005}.
\qed

It is worthy noting that, the 3-dissection of $\overline{R}(\xi_6^2;q)$ plays a key role in \cite{Lovejoy-Osburn-2008}, where
Lovejoy and Osburn elaborated rank differences of overpartitions modulo $3$.
In fact, we have
\begin{align*}
\nonumber\overline{R}(\xi_6^2;q)
&=\sum_{n=0}^\infty(\overline{N}(0,6,n)-\overline{N}(1,6,n)
-\overline{N}(2,6,n)+\overline{N}(3,6,n))q^n\\
&=\sum_{n=0}^\infty(\overline{N}(0,6,n)-\overline{N}(1,6,n)
-\overline{N}(4,6,n)+\overline{N}(3,6,n))q^n.\\
&=\sum_{n=0}^\infty(\overline{N}(0,3,n)-\overline{N}(1,3,n))q^n.
\end{align*}
Thus, we have provided a new proof of \cite[Theorem 1]{Lovejoy-Osburn-2008}
\footnote{Lemma \ref{R(xi_6;q)} reduced $-1$ from \cite[Theorem 1]{Lovejoy-Osburn-2008}, since we assume the convention $\overline{p}(0)=1$.}.

$\overline{R}(\xi_6^3;q)$ is the most complicated part, since double poles occur.
Bringmann and Lovejoy \cite{Bringmann-Lovejoy-2007} pointed out that $\overline{R}(\xi_6^3;q)$ is the holomorphic part of a harmonic weak Maass form of half integral weight.
We here elaborate its property of 3-dissection.
In view of \eqref{xi_6^3}, a straightforward idea is to split the sum into three sums according to the summation index $n$ modulo 3, such as
\begin{align*}
\sum_{n=-\infty}^\infty\frac{(-1)^nq^{n^2+n}}{(1+q^{n})^2}
&=\sum_{n=-\infty}^\infty\frac{(-1)^nq^{9n^2+3n}}{(1+q^{3n})^2}\\
&\indent\indent
-\sum_{n=-\infty}^\infty\frac{(-1)^nq^{9n^2+9n+2}}{(1+q^{3n+1})^2}
+\sum_{n=-\infty}^\infty\frac{(-1)^nq^{9n^2+15n+6}}{(1+q^{3n+2})^2}.
\end{align*}
For the sake of matching the order of $q$ in both numerators and denominators, Lemma \ref{Chan-Thm-3.2} and Theorem \ref{main-result} will
generate identities containing seven or six generalized Lambert series respectively, some of which are redundant.
Therefore, we first make some adjustments in $\overline{R}(\xi_6^3;q)$,
\begin{align*}
\sum_{n=-\infty}^\infty\frac{(-1)^nq^{n^2+n}}{(1+q^n)^2}
&=\sum_{n=-\infty}^\infty\frac{(-1)^nq^{n^2+n}(1-q^n+q^{2n})^2}{(1+q^{3n})^2}
\nonumber\\
&=\sum_{n=-\infty}^\infty\frac{(-1)^nq^{n^2+n}(1-2q^n+3q^{2n}-2q^{3n}+q^{4n})}
{(1+q^{3n})^2}.
\end{align*}
Noting that
\begin{align*}
\sum_{n=-\infty}^\infty\frac{(-1)^nq^{n^2+n}q^{mn}}{(1+q^{3n})^2}
=\sum_{n=-\infty}^\infty\frac{(-1)^nq^{n^2+n}q^{(4-m)n}}{(1+q^{3n})^2},
\end{align*}
we have
\begin{align}\label{4.10}
\nonumber\sum_{n=-\infty}^\infty\frac{(-1)^nq^{n^2+n}}{(1+q^n)^2}
&=2\sum_{n=-\infty}^\infty\frac{(-1)^nq^{n^2+n}}{(1+q^{3n})^2}\\
&\quad-4\sum_{n=-\infty}^\infty\frac{(-1)^nq^{n^2+2n}}{(1+q^{3n})^2}
+3\sum_{n=-\infty}^\infty\frac{(-1)^nq^{n^2+3n}}{(1+q^{3n})^2}.
\end{align}
Then Corollary \ref{cor3.8} will help.
We now give the 3-dissection of $\overline{R}(\xi_6^3;q)$.
\begin{lem}\label{R(xi_6^3;q)}
We have
\begin{align*}
&\overline{R}(\xi_6^3;q)\nonumber\\[5pt]
&=\left(-\frac{2J_{3,18}^4J_6^2J_{9,18}^3}{J_{18}^6}
+\frac{12}{J_{9,18}}\sum_{n=-\infty}^\infty
\frac{(-1)^nq^{9n^2+9n}}{(1+q^{9n})^2}\right)
+q\frac{2J_{3,18}^5J_6^2J_{9,18}^2}{J_{18}^6}\nonumber\\[5pt]
&
+q^2\left(\frac{4J_{3,18}^6J_6^2J_{9,18}}{J_{18}^6}
-\frac{24}{J_{9,18}}\sum_{n=-\infty}^\infty\frac{(-1)^nq^{9n^2+9n}}{(1+q^{9n+3})^2}
+\frac{16}{J_{9,18}}\sum_{n=-\infty}^\infty\frac{(-1)^nq^{9n^2+9n}}{1+q^{9n+3}}\right).
\end{align*}
\end{lem}
\pf
By Corollary \ref{cor3.8} and \eqref{4.10}, we find that
\begin{align*}
&\sum_{n=-\infty}^\infty\frac{(-1)^nq^{n^2+n}}{(1+q^n)^2}
\nonumber\\[5pt]
&\quad
=\frac{J_{3,18}^{12}J_6^6J_{9,18}^4}{4J_{18}^{18}}
+\left(1-\frac{2qJ_{3,18}}{J_{9,18}}\right)
\left(3\sum_{n=-\infty}^\infty\frac{(-1)^nq^{9n^2+9n}}{(1+q^{9n})^2}
-\frac{3J_{3,18}^4J_6^2J_{9,18}^4}{4J_{18}^6}\right)
\nonumber\\[5pt]
&\indent\indent
-q^2\left(1-\frac{2qJ_{3,18}}{J_{9,18}}\right)
\left(6\sum_{n=-\infty}^\infty\frac{(-1)^nq^{9n^2+9n}}{(1+q^{9n+3})^2}
-4\sum_{n=-\infty}^\infty\frac{(-1)^n q^{9n^2+9n}}{1+q^{9n+3}}\right).
\end{align*}
Multiply both sides of $\frac{4(-q)_\infty}{(q)_\infty}$, then we have
\begin{align*}
&\overline{R}(\xi_6^3;q)\\
&=\frac{(-q)_\infty}{(q)_\infty}\frac{J_{3,18}^{12}J_6^6J_{9,18}^4}{J_{18}^{18}}
+\frac{1}{J_{9,18}}\left(12\sum_{n=-\infty}^\infty\frac{(-1)^nq^{9n^2+9n}}{(1+q^{9n})^2}
-\frac{3J_{3,18}^4J_6^2J_{9,18}^4}{J_{18}^6}\right)
\nonumber\\[5pt]
&\indent\indent
-\frac{q^2}{J_{9,18}}
\left(24\sum_{n=-\infty}^\infty\frac{(-1)^nq^{9n^2+9n}}{(1+q^{9n+3})^2}
-16\sum_{n=-\infty}^\infty\frac{(-1)^n q^{9n^2+9n}}{1+q^{9n+3}}\right).
\end{align*}

Compared with Lemma \ref{R(xi_6^3;q)}, it suffices to prove that
\begin{align*}
\frac{(-q)_\infty}{(q)_\infty}\frac{J_{3,18}^{12}J_6^6J_{9,18}^4}{J_{18}^{18}}
=\frac{J_{3,18}^4J_6^2J_{9,18}^3}{J_{18}^6}
+q\frac{2J_{3,18}^5J_6^2J_{9,18}^2}{J_{18}^6}
+q^2\frac{4J_{3,18}^6J_6^2J_{9,18}}{J_{18}^6},
\end{align*}
which is equivalent to \eqref{Hirschhorn-Sellers-2005}.
\qed

Summing up Lemma \ref{R(1;q)}-\ref{R(xi_6^3;q)}, we are now equipped to elaborate 3-dissections of ranks for overpartitions, in view of Lemma \ref{lem4.1}.
Recall that $\overline{r}_s(d)$ is defined in \eqref{r_s(d)}.

\begin{thm}\label{rank-diff-0}
For $d=0$, we have
\begin{align*}
\overline{r}_0(0)
&=\frac{J_{6}^{12}}{6J_{1,6}^8J_{2}^4J_{3,6}}
+\frac{2J_{6}^3J_{3,6}}{3J^2_{1,6}J_{2}}
-\frac{J_{1,6}^4J_2^2J_{3,6}^3}{3J_{6}^6}
+\frac{2}{J_{3,6}}\sum_{n=-\infty}^\infty
\frac{(-1)^nq^{3n^2+3n}}{(1+q^{3n})^2},\\[5pt]
\overline{r}_1(0)
&=\frac{J_{6}^{12}}{6J_{1,6}^8J_{2}^4J_{3,6}}
\indent\indent\indent\indent~
+\frac{J_{1,6}^4J_2^2J_{3,6}^3}{3J_{6}^6}
-\frac{2}{J_{3,6}}\sum_{n=-\infty}^\infty
\frac{(-1)^nq^{3n^2+3n}}{(1+q^{3n})^2},\\[5pt]
\overline{r}_2(0)
&=\frac{J_{6}^{12}}{6J_{1,6}^8J_{2}^4J_{3,6}}
-\frac{J_{6}^3J_{3,6}}{3J^2_{1,6}J_{2}}
-\frac{J_{1,6}^4J_2^2J_{3,6}^3}{3J_{6}^6}
+\frac{2}{J_{3,6}}\sum_{n=-\infty}^\infty
\frac{(-1)^nq^{3n^2+3n}}{(1+q^{3n})^2},\\[5pt]
\overline{r}_3(0)
&=\frac{J_{6}^{12}}{6J_{1,6}^8J_{2}^4J_{3,6}}
\indent\indent\indent\indent~
+\frac{J_{1,6}^4J_2^2J_{3,6}^3}{3J_{6}^6}
-\frac{2}{J_{3,6}}\sum_{n=-\infty}^\infty
\frac{(-1)^nq^{3n^2+3n}}{(1+q^{3n})^2}.
\end{align*}
\end{thm}

\begin{thm}\label{rank-diff-1}
For $d=1$, we have
\begin{align*}
\overline{r}_0(1)
&=\frac{J_{6}^{12}}{3J_{1,6}^7J_{2}^4J_{3,6}^2}
+\frac{4J_{6}^3}{3J_{1,6}J_{2}}
+\frac{J_{1,6}^5J_2^2J_{3,6}^2}{3J_{6}^6},\\[5pt]
\overline{r}_1(1)
&=\frac{J_{6}^{12}}{3J_{1,6}^7J_{2}^4J_{3,6}^2}
\indent\indent\indent\indent~
-\frac{J_{1,6}^5J_2^2J_{3,6}^2}{3J_{6}^6},\\[5pt]
\overline{r}_2(1)
&=\frac{J_{6}^{12}}{3J_{1,6}^7J_{2}^4J_{3,6}^2}
-\frac{2J_{6}^3}{3J_{1,6}J_{2}}
+\frac{J_{1,6}^5J_2^2J_{3,6}^2}{3J_{6}^6},\\[5pt]
\overline{r}_3(1)
&=\frac{J_{6}^{12}}{3J_{1,6}^7J_{2}^4J_{3,6}^2}
\indent\indent\indent\indent~
-\frac{J_{1,6}^5J_2^2J_{3,6}^2}{3J_{6}^6}.
\end{align*}
\end{thm}

The case $d=2$ would be kind of complicated.
Here we use Corollary \ref{cor4.1} to remove the generalized Lambert series with denominator $(1-q^{3n+1})$.

\begin{thm}\label{rank-diff-2}
For $d=2$, we have
\begin{align*}
\overline{r}_0(2)
&=\frac{2J_{6}^{12}}{3J_{1,6}^6J_{2}^4J_{3,6}^3}
-\frac{4J^3_{6}}{3J_{2}J_{3,6}}
+\frac{2J_{1,6}^6J_2^2J_{3,6}}{3J_{6}^6}
\nonumber\\[5pt]
&\quad\quad\quad
-\frac{4}{J_{3,6}}
\sum_{n=-\infty}^\infty
\frac{(-1)^nq^{3n^2+3n}}{(1+q^{3n+1})^2}
+\frac{4}{J_{3,6}}\sum_{n=-\infty}^\infty\frac{(-1)^n q^{3n^2+3n}}{1+q^{3n+1}},\\
\overline{r}_1(2)
&=\frac{2J_{6}^{12}}{3J_{1,6}^6J_{2}^4J_{3,6}^3}
+\frac{2J_{6}^3}{J_{2}J_{3,6}}
-\frac{2J_{1,6}^6J_2^2J_{3,6}}{3J_{6}^6}
\nonumber\\[5pt]
&\quad\quad\quad+\frac{4}{J_{3,6}}
\sum_{n=-\infty}^\infty
\frac{(-1)^nq^{3n^2+3n}}{(1+q^{3n+1})^2}
-\frac{4}{J_{3,6}}\sum_{n=-\infty}^\infty\frac{(-1)^n q^{3n^2+3n}}{1+q^{3n+1}},\\
\overline{r}_2(2)
&=\frac{2J_{6}^{12}}{3J_{1,6}^6J_{2}^4J_{3,6}^3}
+\frac{2J^3_{6}}{3J_{2}J_{3,6}}
+\frac{2J_{1,6}^6J_2^2J_{3,6}}{3J_{6}^6}
\nonumber\\[5pt]
&\quad\quad\quad
-\frac{4}{J_{3,6}}
\sum_{n=-\infty}^\infty
\frac{(-1)^nq^{3n^2+3n}}{(1+q^{3n+1})^2}
+\frac{2}{J_{3,6}}\sum_{n=-\infty}^\infty\frac{(-1)^n q^{3n^2+3n}}{1+q^{3n+1}},\\
\overline{r}_3(2)
&=\frac{2J_{6}^{12}}{3J_{1,6}^6J_{2}^4J_{3,6}^3}
-\frac{4J_{6}^3}{J_2J_{3,6}}
-\frac{2J_{1,6}^6J_2^2J_{3,6}}{3J_{6}^6}
\nonumber\\[5pt]
&\quad\quad\quad
+\frac{4}{J_{3,6}}\sum_{n=-\infty}^\infty
\frac{(-1)^nq^{3n^2+3n}}{(1+q^{3n+1})^2}.
\end{align*}
\end{thm}

Theorem \ref{rank-diff-0}-\ref{rank-diff-2} suggest information of rank sizes for different residues.
Some of the comparisons are quite trivial.
For example, by Theorem \ref{rank-diff-0}, it is easy to derive that, the following inequalities hold for $n\ge 1$:
\begin{align*}
\overline{N}(1,6,3n)=\overline{N}(3,6,3n),\\
\overline{N}(0,6,3n)\ge\overline{N}(2,6,3n).
\end{align*}
Though, some other comparisons take more efforts.
We find that, for fixed $d$, the generating functions of ranks for each residue share a common main term, which is the first term.
After taking differences, the growth rates of the second terms (some have the coefficient $0$) overcome all the others left.
This results in a total ordering relation.
For large integers, this should be able to verify by computing efficient asymptotic formulas of all terms, using standard analytic methods.
While for small ones, this can be verified directly by computer.
Though, this is far from the theme of this article and would take up a dozen pages.
Therefore, we leave it as a conjecture here.

\begin{conj}
For $n\ge 11$, we have
\begin{align*}
\overline{N}(0,6,3n)\ge\overline{N}(1,6,3n)&=\overline{N}(3,6,3n)\ge\overline{N}(2,6,3n),
\\[3pt]
\overline{N}(0,6,3n+1)\ge\overline{N}(1,6,3n+1)&=\overline{N}(3,6,3n+1)\ge\overline{N}(2,6,3n+1),
\\[3pt]
\overline{N}(1,6,3n+2)\ge\overline{N}(2,6,3n+2)
&\ge\overline{N}(0,6,3n+2)\ge\overline{N}(3,6,3n+2).
\end{align*}
\end{conj}

\section{Mock Theta Functions}

Recall that the Appell-Lerch sum is defined as
\begin{equation}
m(x,q,z):=\frac{1}{j(z;q)}\sum_{r=-\infty }^{\infty}
\frac{(-1)^r q^{r\choose 2}z^r}{1-q^{r-1}xz},
\end{equation}
where $x,z \in \mathbb{C}^*$ with neither $z$ nor $xz$ an integral power of $q$.
In \cite{Hickerson-Mortenson-2014}, it points out that the third order mock theta functions $\omega(q)$ and $\rho(q)$ can be expressed in term of $m(x,q,z)$ as follows,
\begin{align}\label{a2}
\omega(q)&=-2q^{-1}m(q,q^{6},q^2)+\frac{J^3_{6}}{J_{2}J_{3,6}},\\[3pt]
\label{a1}
\rho(q)&=q^{-1}m(q,q^{6},-q).
\end{align}

A generalized Lambert series with single poles is essentially an Appell-Lerch sums, so it plays the key role in combining rank differences and mock theta functions.
This section is devoted to proving the relations between the rank differences of overpartitions and mock theta functions, as stated in Theorem \ref{mock}.

First we recall the universal mock theta function $g_2(x,q)$  defined by  Gordon and McIntosh \cite{Gordon-2012}
\[g_2(x,q):=\frac{1}{J_{1,2}}\sum_{n=-\infty}^\infty\frac{(-1)^nq^{n(n+1)}}{1-xq^n}.\]
Hickerson and  Mortenson  \cite{Hickerson-Mortenson-2014} showed that $g_2(x,q)$ and $m(x,q,z)$ have the following relation,
\begin{equation}\label{g2-m}
  g_2(x,q)=-x^{-1}m(x^{-2}q,q^2,x).
\end{equation}
They also introduced the following identities on $m(x,q,z)$:
\begin{align}
m(x,q,z)&=x^{-1}m(x^{-1},q,z^{-1}),\label{mm2} \\[3pt]
m(x,q,z_1)-m(x,q,z_0)&=\frac{z_0J_1^3j(z_1/z_0;q)j(xz_0z_1;q)}
{j(z_0;q)j(z_1;q)j(xz_0;q)j(xz_1;q)}.\label{mm}
\end{align}

We are now in a position to give a proof of Theorem  \ref{mock}. \vskip 0.2cm

\noindent{\it Proof of Theorem \ref{mock}.}
From Theorem \ref{rank-diff-2}, we have
\begin{align}\label{mock-1}
\overline{r}_0(2)+\overline{r}_3(2)
&=\frac{4}{J_{3,6}}\sum_{n=-\infty}^\infty\frac{(-1)^n q^{3n^2+3n}}{1+q^{3n+1}}
+\frac{4J_{6}^{12}}{3J_{1,6}^6J_{2}^4J_{3,6}^3}
-\frac{16J^3_{6}}{3J_{2}J_{3,6}}.
\end{align}
Replacing $q$ by $q^3$ in \eqref{g2-m} and setting $x=-q$, we have
\begin{equation}
g_2(-q,q^3)=q^{-1}m(q,q^6,-q),
\end{equation}
and by \eqref{a1}, we deduce that
\begin{equation}\label{m-j-1}
\rho(q)=g_2(-q,q^3).
\end{equation}
Together with the identity in \cite[p.63]{Watson-1936}
\begin{equation}\label{m-j}
\omega(q)+2\rho(q)=\frac{3J^3_{6}}{J_{2}J_{3,6}},
\end{equation}
we find that \eqref{mock-1} can be transformed as follows:
\begin{align*}\nonumber
\overline{r}_0(2)+\overline{r}_3(2)
&=4\rho(q)-\frac{16}{9}(\omega(q)+2\rho(q))
+\frac{4J_{6}^{12}}{3J_{1,6}^6J_{2}^4J_{3,6}^3}
\\[3pt]
&=\frac{4}{9}\rho(q)-\frac{16}{9}\omega(q)+\frac{4J_{6}^{12}}{3J_{1,6}^6J_{2}^4J_{3,6}^3}.
\end{align*}

Similarly, we have
\begin{align*}
\overline{r}_1(2)-\overline{r}_3(2)
&=-\frac{4}{J_{3,6}}\sum_{n=-\infty}^\infty\frac{(-1)^n q^{3n^2+3n}}{1+q^{3n+1}}
+\frac{6J^3_{6}}{3J_{2}J_{3,6}}
\\[3pt]
&=2\omega(q),
\end{align*}
and
\begin{align*}
\overline{r}_2(2)+\overline{r}_3(2)
&=\frac{2}{J_{3,6}}\sum_{n=-\infty}^\infty\frac{(-1)^n q^{3n^2+3n}}{1+q^{3n+1}}
-\frac{10J^3_{6}}{3J_{2}J_{3,6}}+\frac{4J_{6}^{12}}{3J_{1,6}^6J_{2}^4J_{3,6}^3}
\\[3pt]
&=2\rho(q)-\frac{10}{9}(\omega(q)+2\rho(q))+\frac{4J_{6}^{12}}{3J_{1,6}^6J_{2}^4J_{3,6}^3}
\\[3pt]
&=-\frac{2}{9}\rho(q)-\frac{10}{9}\omega(q)+\frac{4J_{6}^{12}}{3J_{1,6}^6J_{2}^4J_{3,6}^3}.
\end{align*}

Thus we complete the proof of Theorem \ref{mock}.\qed

%
%


\begin{thebibliography}{99} \small
\setlength{\itemsep}{-.8mm}


\bibitem{Andrews-Hickerson-1991}
G. E. Andrews and D. Hickerson, Ramanujan's ``lost'' notebook VII: The sixth order mock theta functions,
Adv. Math. 89 (1991) 60--105.

\bibitem{Andrews-Lewis-2000}
G. E. Andrews and R. Lewis, The ranks and cranks of partitions moduli $2,3$, and $4$, J. Number Theory 85 (2000) 74--84.

\bibitem{Andrews-Lewis-Liu-2001}
G. E. Andrews, R. Lewis and Z.-G. Liu, An identity relating a theta function to a sum of Lambert series, Bull. London Math. Soc. 33 (1) (2001) 25--31.


\bibitem{Atkin-Swinnerton-Dyer-1954}
A. O. L. Atkin and P. Swinnerton-Dyer, Some properties of partitions, Proc. Lond. Math. Soc. 66 (1954) 84--106.

\bibitem{Bringmann-Lovejoy-2007}
K. Bringmann and J. Lovejoy, Dyson's rank, overpartitions, and weak Maass forms, Int. Math. Res. Not. 19 (2007) Art. ID rnm063 34 pp.

\bibitem{Chan-2005}
S. H. Chan, Generalized Lambert series identities, Proc. London Math. Soc. 91 (3) (2005) 598--622.


\bibitem{Garvan-1988}
F. G. Garvan, New combinatorial interpretations of Ramanujan's partition congruences mod 5, 7 and 11, Trans. Amer. Math. Soc. 305 (1988) 47--77.

\bibitem{Garvan-1990}
F. G. Garvan, The crank of partitions mod 8,9 and 10, Trans. Amer. Math. Soc. 322 (1) (1990) 79--94.

\bibitem{Gordon-2012}
B. Gordon and R. J. McIntosh, A survey of classical mock theta functions, In:  K. Alladi  and F. G. Garvan,
Partitions, $q$-Series, and Modular Forms, Developmental Mathematics  23 (2012) 95--144.

\bibitem{Hickerson-Mortenson-2014}
D. Hickerson and E. Mortenson, Hecke-type double sums, Appell-Lerch sums, and mock theta functions, \textrm{I}, Proc. Lond. Math. Soc. 109 (3) (2014) 382--422.

\bibitem{Hirschhorn-Sellers-2005}
M. D. Hirschhorn and J. A. Sellers, Arithmetic relations for overpartitions, J. Combin. Math. Combin. Comput. 53 (2005) 65--73.

\bibitem{Jennings-Shaffer-2016}
C. Jennings-Shaffer, Overpartition rank differences modulo $7$ by Maass forms, J. Number Theory 163 (2016) 331--358.

\bibitem{Ji-Zhang-Zhao-2017}
K. Q. Ji, H. W. J. Zhang and A. X. H. Zhao, Ranks of overpartitions modulo $6$ and $10$, J. Number Theory 184 (2018) 235--269.




\bibitem{Lovejoy-2005}
J. Lovejoy, Rank and conjugation for the Frobenius representation of an overpartition, Ann. Combin. 9 (3) (2005) 321--334.

\bibitem{Lovejoy-Osburn-2008}
J. Lovejoy and R. Osburn, Rank differences for overpartitions, Quart. J. Math. 59 (2) (2008) 257--273.


\bibitem{Mao-2013}
R. R. Mao, Ranks of partitions modulo $10$, J. Number Theory 133 (2013) 3678--3702.

\bibitem{Mao-2015}
R. R. Mao, The $M_2$-rank of partitions without repeated odd parts modulo $6$ and $10$, Ramanujan J. 37 (2015) 391--419.


\bibitem{Sears-1951}
D. B. Sears, On the transformation theory of hypergeometric functions and cognate trigonometric series, Proc. London Math. Soc. 53 (2) (1951) 138--157.

\bibitem{Watson-1936}
G. N. Watson, The final problem: an account of the mock theta functions, J. London Math. Soc. 11 (1936) 55--80.


\end{thebibliography}
\end{document}